\documentclass[11pt]{article}
\usepackage{amsfonts}
\usepackage{amssymb}
\usepackage{amsmath}
\usepackage[english]{babel}
\usepackage{amssymb}
\usepackage{graphicx}
\usepackage{graphics}
\usepackage{amsthm}
\newcommand{\ZZ}{\mathbb{Z}}

\newcommand{\esp}{\thinspace}
\theoremstyle{definition}

\input epsf
\addto\captionsspanish{}

\begin{document}

\date{}
\author{Andr\'es Navas}

\title{Sur les rapprochements par conjugaison en dimension 1 et classe $C^1$}
\maketitle

\vspace{-0.2cm}

Toute action d'un groupe d\'enombrable par hom\'eomorphismes d'une vari\'et\'e
unidimensionnelle (s\'epar\'ee) est topologiquement conjugu\'ee \`a une action par des
hom\'eomorphismes bilipschitziens. Ce r\'esultat \`a l'air innocent a \'et\'e \'etabli
dans \cite{DKN-acta} via une m\'ethode probabiliste pour des vari\'et\'es compactes,
mais une preuve plus simple et g\'en\'erale a \'et\'e donn\'ee par B.~Deroin dans
\cite{deroin}. Il est important de signaler que ceci est loin d'\^etre valable
en dimension sup\'erieure m\^eme pour des actions de $\mathbb{Z}$,
d'apr\`es notamment \cite{harrison}.

Les r\'esultats de cette Note sont inspir\'es (entre autres) par le fait ci-dessus
ainsi que par les questions suivantes :

\noindent {\em -- Sous quelles conditions une action donn\'ee peut-elle 
\^etre conjugu\'ee en une action dont les g\'en\'era- teurs deviennent 
aussi (Lipschitz ou $C^1$) proches de translations que l'on veut ? 

\noindent -- 
Dans le cas o\`u de telles conjugaisons existent, 
peut-on les relier par un chemin continu de conjugu\'es ?} 

Voici un premier r\'esultat dans le contexte lipschitzien qui nous 
sert de motivation pour la suite. Pour simplifier, nous d\'esign\'erons
par $X$ soit le cercle soit l'intervalle ferm\'e, et par la suite nous ne 
consid\'ererons que des hom\'eomorphismes de $X$ qui respectent l'orientation.

\vspace{0.38cm}

\noindent{\bf Th\'eor\`eme A.} Si $\Gamma$ est un groupe de type fini et \`a
croissance sous-exponentielle d'hom\'eomorphismes de $X$, alors pour tout
$\varepsilon > 0$ il existe des conjugu\'es topologiques de $\Gamma$ pour
lesquels les g\'en\'erateurs (et leurs inverses) sont des hom\'eomorphismes
lipschitziens \`a des contantes de Lipschitz $\leq e^{\varepsilon}$.

\vspace{0.38cm}

Dans le cadre des diff\'eomorphismes de classe $C^1$, nous ne savons pas
traiter en g\'en\'eral les actions de groupes \`a croissance sous-exponentielle.
Cependant, nous pouvons donner une r\'eponse affirmative \`a nos questions pour
les groupes nilpotents. Signalons d'une part que des constructions diff\'erentes
d'actions de groupes (sans torsion et) nilpotents par diff\'eomorphismes de
classe $C^1$ de l'intervalle (ainsi que des r\'esultats de rigidit\'e en classe
$C^{1+\tau}$) sont donn\'ees dans \cite{CJN,FF,jorquera,JAM}. D'autre part, un
exemple d'un groupe de diff\'eomorphismes \`a croissance sous-exponentielle et
non virtuellement nilpotent est donn\'e dans \cite{gafa}. 

\vspace{0.38cm}

\noindent{\bf Th\'eor\`eme B.} Si $\Gamma$ est un groupe nilpotent et de type 
fini de diff\'eomorphismes de classe $C^1$ de $X$, alors pour tout $\varepsilon > 0$ il
existe des conjugu\'es topologiques de $\Gamma$ pour lesquels les g\'en\'erateurs
(et leurs inverses) sont encore des diff\'eomorphismes de classe $C^1$ de $X$
mais \`a d\'eriv\'ee $\leq e^{\varepsilon}$ partout.

\vspace{0.38cm}

Rappelons que pour une action sur $X$, l'orbite d'une paire de points $x < y$
est dite {\em de type ressort} s'il existe des \'el\'ements $f,g$ tels que
\begin{equation}\label{resilient}
x < f(x) < f(y) < g(x) < g(y) < y.
\end{equation}
L'existence de telles orbites est une obstruction pour rapprocher (au sens de Lipschitz)
des actions par des translations. En effet, les relations (\ref{resilient}) sont
stables par conjugaison topologique, et elles entra\^{\i}nent \'evidemment que
l'un des \'el\'ements doit contracter d'un facteur $< 1/2$ une partie de
l'intervalle correspondant. Encore plus, ces orbites donnent lieu \`a de
l'entropie positive pour l'action \cite{GLW}. D'un cot\'e alg\'ebrique,
les orbites de type ressort entra\^{\i}nent l'existence de semigroupes
libres (par une application directe du lemme du ping-pong de Klein dans
sa version positive), donc une croissance exponentielle pour le groupe.
Ceci rend naturelle la question suivante :

\vspace{0.2cm}

\noindent{\bf Question.} L'absence d'orbites de type ressort pour une action
par diff\'eomorphismes de classe $C^1$ entra\^{\i}ne-t-elle l'existence
de conjugaisons topologiques de telle sorte que les g\'en\'erateurs
deviennent aussi Lipschitz ($C^1$ ?) proches de translations que l'on veut ?

\vspace{0.2cm}

Cette question se pose plus naturellement dans le contexte des pseudo-groupes
d'hom\'eomorphis- mes (donc, pour des feuilletages de codimension 1). Dans
ce cadre, une r\'eponse par l'affirmatif donnerait une preuve
alternative du r\'esultat de \cite{hurder} qui \'etablit la nullit\'e de
l'entropie g\'eom\'etrique \cite{GLW} de tout feuilletage de codimension 1,
transversalement $C^1$ et sans feuille ressort. Signalons en passant qu'il
existe des actions lipschitziennes \`a entropie positive et sans orbite
ressort \cite{GLW}, mais ces actions ne sont pas $C^1$ lissables \cite{CC}.

\vspace{0.1cm}

Pour conclure, nous consid\'erons une version am\'elior\'ee du Th\'eor\`eme B, 
o\`u les conjugu\'es forment un chemin continu pour la topologie $C^1$. Ceci 
est \'etroitement li\'e \`a des probl\`emes de d\'eformation de feuilletages de 
codimension 1 (voir \cite{eynard2}).\footnote{Une autre motivation vient d'une 
vieille question de H. Rosenberg \`a propos de l'existence d'actions structurellement 
stables de $\mathbb{Z}^2$ par diff\'eomorphismes du cercle. \`A notre connaissance, 
cette question reste encore ouverte.} Dans cet esprit, dans 
\cite{eynard}, H.~Eynard s'int\'eresse aux r\'epr\'esentations de $\ZZ^d$ dans
$\mathrm{Diff}_+([0,1])$ et d\'emontre la $C^1$-connexit\'e par
arcs de l'espace des r\'epr\'esentations par diff\'eomorphismes
de classe $C^2$ (voir la pr\'epublication r\'ecente \cite{BE}
pour la connexit\'e par arcs en classe $C^{\infty}$). Le
th\'eor\`eme plus bas \'etend ce r\'esultat en ce qui concerne
le groupe qui agit, la r\'egularit\'e des diff\'eomorphismes
concern\'es et la vari\'et\'e unidimensionnelle sous-jacente.

\vspace{0.38cm}

\noindent{\bf Th\'eor\`eme C.}  L'espace des actions de tout groupe nilpotent de type 
fini par diff\'eomorphismes de classe $C^{1}$ de l'intervalle ferm\'e est connexe par arcs. 
Ceci reste valable pour le cas du cercle pour des actions de $\mathbb{Z}^d$ ; plus 
g\'en\'eralement, toute action d'un groupe nilpotent de type fini par diff\'eomorphismes 
de classe $C^{1}$ du cercle est dans la composante connexe d'une action qui transite 
par un morphisme vers un groupe fini de rotations.

\vspace{0.38cm}

Au del\`a des groupes nilpotents, nous ne connaissons pas d'autres groupes {\em moyennables} 
pours lesquels ce th\'eor\`eme reste valable. Dans cette direction, deux exemples int\'eressants 
\`a traiter pour commencer ce sont le groupe de Grigorchuk-Mach\`{\i} \cite{gafa} 
et le groupe de Baumslag-Solitar $BS(1,2)$.  Notons cependant que dans la preuve 
du Th\'eor\`eme C, l'arc qui joint deux r\'epr\'esentations transite par une 
r\'epr\'esentation par des translations. Pour joindre une r\'epr\'esentation donn\'ee 
\`a une r\'epr\'esentation par des translations, nous construisons un chemin
explicite form\'e par des conjugu\'es topologiques de la r\'epr\'esentation
originelle. Or, pour le cas du groupe $BS(1,2)$, les r\'esultats de 
\cite{CC,nancy,rivas} plus la discussion autour des orbites de type
ressort plus haut impliquent que l'action triviale ne peut \^etre
rapproch\'ee (en topologie $C^1$) par des conjugu\'es topologiques 
d'aucune action fid\`ele par diff\'eomorphismes de classe $C^1$.

Bien \'evidemment, le th\'eor\`eme est valable aussi pour des actions du groupe libre. En 
effet, dans ce contexte, un argument simple de transversalit\'e montre que l'espace des 
actions {\em fid\`eles} est connexe par arcs. Or, ces chemins ne peuvent pas toujours 
venir de conjugu\'es topologiques des actions donn\'ees. Par exemple, une action de 
type Schottky sur le cercle ne peut jamais rapprocher --m\^eme contin\^ument-- 
aucune action par des rotations via des conjugaisons topologiques.


Pour conclure cette Introduction, signalons qu'\'etant donn\'ee la m\'ethode
de d\'emonstration du Th\'eor\`eme C d\'ecrite plus haut, nous n'obtenons
aucune information de connexit\'e {\em locale} par arcs pour l'espace
de r\'epr\'esentations. Ce probl\`eme reste largement ouvert.

\vspace{0.3cm}

\noindent{\bf I. Sur les actions de groupes \`a croissance sous-exponentielle.}
Pour la preuve du Th\'eor\`eme A, nous suivons la m\'ethode de \cite{deroin}.
Fixons $\lambda := e^{-\varepsilon} < 1$. Puisque $\Gamma$ a une croissance
sous-exponentielle, \'etant fix\'e un syst\`eme fini de g\'en\'erateurs $\mathcal{G}$,
il existe $C = C_{\varepsilon,\mathcal{G}}$ tel que le cardinal de la boule
$B(n)$ de rayon $n$ correspondante est $\leq C([ \lambda + 1] / 2 \lambda)^n$.
Consid\'erons la mesure $\mu$ sur $X$ d\'efinie par
$$\mu := \sum_{f \in \Gamma} \lambda^{\ell(f)} f_* (Leb),$$
o\`u $\ell(f)$ d\'esigne la longueur de l'\'el\'ement $f$ par rapport au syst\`eme de
g\'en\'erateurs choisi et $Leb$ d\'esigne la mesure de Lebesgue sur $X$. Nous
affirmons que $\mu$ a une masse totale finie. En effet, si l'on d\'esigne par
$S(n)$ la sph\`ere de rayon $n$ dans $\Gamma$, alors
$$\mu (X) = \sum_{n \geq 0} \lambda^n \big| S(n) \big|
\leq \sum_{n \geq 0} \lambda^n \big| B(n) \big|
\leq C \sum_{n \geq 0} \left( \frac{\lambda + 1}{2} \right)^n
= \frac{2C}{1 - \lambda}
< \infty.$$
De plus, puisque pour tout g\'en\'erateur $g$ et tout $f \in \Gamma$ l'in\'egalit\'e  
\hspace{0.03cm} $\big| \ell(gf) - \ell(f) \big| \leq 1$ \hspace{0.03cm} a lieu, nous avons
\begin{equation}\label{quasi}
g_* (\mu) = \sum_{f \in \Gamma} \lambda^{\ell(f)} (gf)_* (Leb) \leq
\frac{1}{\lambda} \sum_{f \in \Gamma} \lambda^{\ell(gf)} (gf)_* (Leb)
= \frac{\mu}{\lambda}.\end{equation}
La mesure $\mu$ est de masse finie, \`a support totale et sans at\^ome. Elle
est donc \'equivalente par conjugaison topologique \`a la mesure de Lebesgue
\`a un facteur pr\`es. Apr\`es un changement de coordon\'ees envoyant
$\frac{\mu}{\mu(X)}$ sur $Leb$, la relation (\ref{quasi}) devient,
pour tout intervalle $I \subset X$,
$$|g^{-1}(I)| = g_* (Leb)(I) \leq \frac{Leb(I)}{\lambda} = \frac{|I|}{\lambda}.$$
Ceci entra\^{\i}ne que dans ces nouvelles coordonn\'ees, $g^{-1}$ est lipschitzien
de constante $\leq 1 / \lambda$.

\vspace{0.1cm}

\noindent{\bf Remarque.} Nous ignorons si dans le cas d'une action par
hom\'eomorphismes lipschitziens, le conjugant ({\em i.e.} le changement de
coordonn\'ees) ci-dessus peut toujours \^etre pris lipschitzien (voir \cite{weiss}
pour un r\'esultat qui, d'apr\`es le \S II. plus bas, pointe dans une direction
plut\^ot n\'egative). En ce qui concerne les Th\'eor\`emes B et C, la n\'ecessit\'e
de consid\'erer des conjugaisons topologiques vient des points fixes hyperboliques,
dont on ne peut pas se d\'ebarrasser par des conjugaisons lisses.

\vspace{0.3cm}

\noindent{\bf II. Des rapprochements par des conjugaisons lisses via
des \'equations cohomologiques.} Un diff\'eomorphisme $f$ du cercle (resp.
de l'intervalle) est $C^1$-proche de la rotation d'angle $\rho(f)$ (resp. de
l'identit\'e) si et seulement si sa d\'eriv\'ee $Df$ est partout proche de 1. Donc,
pour obtenir des rapprochements via des conjugaisons par diff\'eomorphismes
soit \`a une rotation soit \`a l'identit\'e, nous devons chercher $\varphi$
de telle sorte que
$$\big| \log D(\varphi \circ f \circ \varphi^{-1}) \big| =
\big| \log(D\varphi) \circ (f \circ \varphi^{-1}) +  \log (Df) \circ (\varphi^{-1})
- \log(D\varphi) \circ (\varphi^{-1}) \big|$$
soit partout petit. En d'autres termes, nous devons chercher des solutions
raproch\'ees $u = \log(D\varphi)$ de {\em l'\'equation cohomologique}
\begin{equation}
u - u \circ f = \log (D f).
\label{equa}\end{equation}
Il se trouve que d'apr\`es \cite{ollagnier-pinchon}, ces solutions raproch\'ees
existent non seulement pour un diff\'eomorphisme mais aussi pour des groupes
presque nilpotents pourvu que certains exposants de Lyapunov associ\'es soient tous
nuls.\footnote{Pour le cas d'une seule application, 
ceci r\'esulte d'une application directe du th\'eor\`eme de Hahn-Banach lorsque la 
moyenne de $\log(D f)$ est nulle par rapport \`a toute probabilit\'e invariante.} 
Dans notre contexte, cela appara\^{i}t explicitement dans le cas (ii) du lemme 
ci-dessous.

\vspace{0.2cm}

\noindent{\bf Lemme (d'existence de solutions rapproch\'ees).} {\em Soit $\Gamma$ 
un sous-groupe presque nilpotent et de type fini de $\mathrm{Diff}_+^1(X)$
engendr\'e par une partie finie $\mathcal{G}$. Alors :}

\noindent (i) {\em soit $\Gamma$ admet des orbites finies,}

\noindent (ii) {\em soit $\Gamma$ n'admet pas de telles orbites mais
il est topologiquement semi-conjugu\'e \`a un groupe de rotations. 
De plus, 
pour tout $\varepsilon > 0$ il
existe une fonction continue $u$ telle que pour tout $f \in \mathcal{G}$
l'in\'egalit\'e \esp $|u - u \circ f - \log (Df)| \leq \varepsilon$ \esp
est partout satisfaite.}

\vspace{0.14cm}

\noindent{\bf Preuve.} Pour des actions sur l'intervalle, on est bien s\^ur toujours dans
le cas (i). Pour des actions sur le cercle, la dichotomie entre (i) et (ii) s'applique plus
g\'en\'eralement aux actions de groupes moyennables : voir \cite[Lemma 4.1.2]{navas-book}.
Passons maintenant \`a l'existence de solutions rapproch\'ees de l'\'equation 
cohomologique dans le cas (ii). Pour cela, notons qu'une application directe de 
\cite[Th\'eor\`eme 1]{ollagnier-pinchon} au cocycle $f \mapsto \log (D f)$ montre 
que pour assurer l'existence de ces solutions, nous devons montrer que pour tout 
\'el\'ement $g \in \Gamma$ et toute mesure de probabilit\'e $\mu$ invariante par 
$\Gamma$, nous avons
$$\int_X \log Dg (x) \hspace{0.01cm} d \mu(x) = 0.$$
Nous affirmons que ceci est toujours valable dans le cas (ii). En effet, si $g$ a un nombre de 
rotation irrationnel, alors la moyenne du logarithme de sa d\'eriv\'ee par rapport \`a l'unique 
mesure de probabilit\'e invariante est nulle (voir \cite[Proposition I.I, Chapitre VI]{herman}).  
Pour le cas de nombre de rotation rationnel, notons d'abord que puisque $\Gamma$ 
est semi-conjugu\'e \`a un groupe infini de rotations, il existe une 
unique mesure de probabilit\'e sur $\mathrm{S}^1$ qui est support\'ee sur l'unique 
ensemble non vide compact invariant et minimal de l'action (cet ensemble $K$ soit 
il co\"{\i}ncide avec tout le cercle, soit il est hom\'eomorphe \`a l'ensemble de Cantor). 
De plus, tout \'el\'ement ayant des points fixes doit fixer chaque point de $K$. En 
particulier, si $g \in \Gamma$ a un nombre de rotation rationnel, alors pour un certain 
$N \geq 1$ on a que $g^N$ fixe tous les points de $K$. Par suite, la d\'eriv\'ee de 
$g^N$ est \'egale \`a 1 partout sur $K = \mathrm{supp}(\mu)$, donc
$$0 = \int_{\mathrm{S}^1} \log D g^N (x) \hspace{0.01cm} d\mu(x) =
\sum_{i=0}^{N-1} \int_{\mathrm{S}^1} \log Dg (g^i(x)) \hspace{0.01cm} d\mu(x)
= N \int_{\mathrm{S}^1} \log Dg (x) \hspace{0.01cm} d \mu(x).$$
Par cons\'equent, \hspace{0.04cm} 
$\int_{\mathrm{S}^1} \log Dg (x) \hspace{0.01cm} d \mu(x) = 0,$
\hspace{0.04cm} tel que nous le voulions. $\hfill\square$

\vspace{0.3cm}

\`A l'aide du lemme pr\'ec\'edent, nous pouvons traiter le probl\`eme 
du rapprochement par conjugaison dans le cas (ii) plus haut (l'autre cas
sera trait\'e plus tard). Pour chaque $n \geq 1$ et chaque $s \in [0,1]$, posons
\hspace{0.1cm} 
$$v_{n+s} := (1-s) u_n + s u_{n+1} + C_{n+s},$$ 
o\`u $u_n$ est une fonction partout v\'erifiant
$$\big| u_n - u_n \circ f - \log (Df) \big| \leq \frac{1}{n}$$
pour tout $f \in \mathcal{G}$ et $C_{n+s}$ est l'unique constante qui satisfait
$$\int_X \exp \big( (1-s)u_n + su_{n+1} \big) = \exp(-C_{n+s}).$$
On v\'erifie alors que $x \to \int_0^x \exp(v_{n+s})$ d\'efinit un
diff\'eomorphisme $\varphi_{n+s}$ de $X$ qui varie contin\^ument par rapport
au param\^etre. De plus, $\varphi_1 = Id$, et en renversant les calculs
pr\'ec\'edents, on constante ais\'ement que \hspace{0.03cm}
$\big|\log D(\varphi_{n+s} \circ f \circ \varphi_{n+s}^{-1}) \big|$
\hspace{0.03cm} converge uniform\'ement vers z\'ero lorsque $n$
tends vers l'infini pour tout $f \in \mathcal{G}$, tel qu'on le souhaitait.

Nous venons donc de montrer que toute action v\'erifiant la condition (ii)
contient une r\'epr\'esenta- tion (non n\'ecessairement fid\`ele !) par des
rotations dans son adh\'erence par conjugaisons $C^1$. De plus, cette
r\'epr\'esentation est aboutie par un chemin continu de conjugu\'es. 
Il nous reste donc \`a joindre deux r\'epr\'esentations quelconques 
par des rotations. Or, pour le cas des actions de $\mathbb{Z}^d$, 
cela se fait tout simpl\'ement en faisant bouger les angles 
associ\'es aux g\'en\'erateurs. Pour le cas d'un groupe 
nilpotent quelconque, les angles associ\'es aux  
g\'en\'erateurs d'ordre infini peuvent encore \^etre ram\'en\'es 
\`a z\'ero de mani\`ere continue, ce qui montre que la repr\'esentation 
est dans la composante connexe d'une action qui transite par un 
morphisme vers un groupe fini de rotations. 
Ceci \'etablit donc les Th\'eor\`emes B et C dans le cas (ii).

\vspace{0.25cm}

\noindent{\bf Remarque.} Les calculs pr\'ec\'edents ne peuvent pas \^etre 
renvers\'es de fa\c{c}on \`a montrer par exemple que toute action libre 
par diff\'eomorphismes du cercle est contenue dans l'adh\'erence par
conjugaisons de la r\'epr\'esentation par des rotations correspondante.
En effet, celui-ci est un probl\`eme qui ne se mod\`ele pas par des
\'equations cohomologiques. D'ailleurs, nous ignorons si cela est 
toujours vrai. Signalons cependant que c'est le cas pour des actions de
$\mathbb{Z}$, d'apr\`es un r\'esultat r\'ecent de C.~Bonatti et N.~Guelman
\cite{bona}. Pour le cas de l'intervalle, un r\'esultat analogue a
\'et\'e aussi r\'ecemment prouv\'e par \'E.~Farinelli \cite{fari}.

\vspace{0.38cm}

\noindent{\bf III. \`A propos des solutions rapproch\'ees.} La
preuve pr\'ec\'edente est un peu obscure car elle fait appel \`a
\cite{ollagnier-pinchon}. Pour la commodit\'e du lecteur, nous rendons 
explicite l'argument pour $\Gamma \sim \mathbb{Z}^d$, ce qui
nous permettra de mieux expliquer notre m\'ethode dans le
cas des orbites finies. De plus, pour la commodit\'e du lecteur, 
une nouvelle preuve (d'une version plus g\'en\'erale mais sous 
une hypoth\`ese l\'eg\`erement plus forte) du r\'esultat 
de \cite{ollagnier-pinchon} sera donn\'ee plus bas. 

Pour chaque $n \geq 1$, on consid\`ere la boule positive $B_+(n)$
de rayon $n$ dans $\mathbb{Z}^d$ par rapport au syst\`eme canonique
des g\'en\'erateurs $\mathcal{G} := \{f_1,\ldots,f_d\}$, c'est-\`a-dire
l'ensemble 
$$B_+(n) := \{ f_1^{k_1} \cdots f_d^{k_d} : 0 \leq k_i < n,
\hspace{0.1cm} 1 \leq i \leq d \}.$$ 
On pose
\begin{equation}
u_n (x) := \frac{1}{|B_+(n)|} \sum_{f \in B_+(n)} \log Df (x).
\label{definicion}
\end{equation}
Nous avons
\begin{eqnarray*}
u_n(f_i(x))
&=& \frac{1}{|B_+(n)|} \sum_{f \in B_+(n)} \log Df (f_i(x))\\
&=& \frac{1}{|B_+(n)|} \sum_{f \in B_+(n)} \left[ \log D(f \circ f_i)(x) - \log Df_i(x) \right]\\
&=& - \log Df_i (x) + \frac{1}{|B_+(n)|} \sum_{f \in B_+(n) f_i} \log Df(x),
\end{eqnarray*}
o\`u \hspace{0.02cm} $B_+(n)f_i := \{gf_i: g \in B_+(n)\}$. Donc, la valeur de l'expression
\begin{equation}
\big| u_n(x) - u_n(f_i(x)) - \log Df_i (x) \big|
\label{cercana}
\end{equation}
est inf\'erieure ou \'egale \`a
\vspace{-1cm}
\begin{multline*}
\frac{1}{n^d} \left| \sum_{f \in B_+(n)} \log Df(x)
- \sum_{f \in B_+(n) f_i} \log Df (x) \right| \\
= \hspace{0.15cm} \frac{1}{n^d} \left| \sum_{0 \leq m_j < n, j \neq i}
 \! \big[ \log D(f_1^{m_1} \cdots f_i^0 \cdots f_d^{m_d}) (x)
- \log D(f_1^{m_1} \cdots f_i^n \cdots f_d^{m_d}) (x) \big] \right| \\
= \hspace{0.15cm} \frac{1}{n^d} \left| \sum_{0 \leq m_j < n, j \neq i}
- \log D(f_i^n) (f_1^{m_1} \cdots f_{i-1}^{m_{i-1}}
f_{i+1}^{m_{i+1}} \cdots f_d^{m_d} (x)) \right| \\
= \hspace{0.15cm} \frac{1}{n^d} \left| \sum_{0 \leq m_j < n}
- \log Df_i (f_1^{m_1} \cdots f_{i}^{m_{i}} \cdots f_d^{m_d} (x)) \right|\\
= \hspace{0.15cm}
\left| \int_{X} \log Df_i (y) \hspace{0.15cm} d\mu_{n,x} (y) \right|,\\
\vspace{0.3cm}
\hspace{-9.7cm} \mbox{o\`u } \mu_{n,x} \mbox{ d\'esigne la mesure de probabilit\'e}\\
\end{multline*}
\vspace{-1.2cm}
$$\frac{1}{|B_+(n)|} \sum_{f \in B_+(n)} \delta_{f(x)}$$
(notons que la commutativit\'e --en non la simple moyennabilit\'e-- a \'et\'e
utilis\'ee pour obtenir les deuxi\`eme et troisi\`eme \'egalit\'es plus haut).
Or, si $\mu$ est un point d'accumulation d'une suite de mesures $\mu_{n,x_n}$,
alors $\mu$ est invariante par l'action de $\Gamma$, et l'int\'egrale
en consid\'eration pour $x_n : = x$ converge vers \hspace{0.02cm}
$\int_{X} \log (Df_i) \hspace{0.15cm} d \mu.$ \hspace{0.1cm} Donc, si ces moyennes 
sont suppos\'es d'\^etre toutes nulles, alors la valeur de (\ref{cercana}) 
tend vers z\'ero. En d'autres termes, les fonctions $(u_n)$ forment 
une suite explicite de solutions rapproch\'ees de notre \'equation 
cohomologique, ce qui nous permet de conjuguer l'action donn\'ee en une 
action pr\^oche d'une action par des translations et ainsi conclure la preuve.  

\vspace{0.22cm}

L'argument plus haut montre qu'une condition suffisante (at d'ailleurs n\'ecessaire) pour
l'existen- ce de solutions rapproch\'ees pour notre \'equation cohomologique
est la nullit\'e de l'int\'egrale de chaque fonction
\hspace{0.05cm} $\log (Df_i)$ \hspace{0.05cm} par rapport
\`a toute mesure de probabilit\'e invariante. De plus, ce calcul montre 
que si cette condition n'est pas satisfaite mais ces int\'egrales sont petites
en valeur absolue, alors on peut trouver des solutions \`a des erreurs petits. 
Ceci est rendu explicite dans le lemme ci-dessous,  
qui \'etend \cite[Th\'eor\`eme 1]{ollagnier-pinchon} 
sous une hypoth\`ese l\'eg\`erement plus forte 
(dans \cite{ollagnier-pinchon}, l'hypoth\`ese porte seulement sur les 
probabilit\'es invariantes par $\Gamma$). Cependant, cette version 
sera suffisante pour nos besoins, et elle a l'avantage d'admettre une preuve 
relativement simple avec des outils modernes.\footnote{En suivant la 
m\'ethode de \cite{ollagnier-pinchon}, on peut montrer l'enonc\'e sous 
l'hypoth\`ese plus faible ne comportant que les mesures invariantes 
par $\Gamma$.} 

\vspace{0.286cm}

\noindent{\bf Lemme (d'existence de solutions \'erron\'ees).}  
{\em Soit $\Gamma$ un groupe nilpotent engendr\'e par une
partie finie $\mathcal{G}$. Il existe une partie g\'en\'eratrice finie $\mathcal{G}'$ 
contenant $\mathcal{G}$ telle que pour tout $\varepsilon > 0$ il existe
$\delta > 0$ qui satisfait la propri\'et\'e suivante : si
$\Phi \!: \Gamma \to \mathrm{Diff}_+^1(X)$ est une r\'epr\'esentation
telle que pour tout $f \in \mathcal{G}'$ la valeur absolue de la moyenne de 
$\log D\Phi(f)$ par rapport \`a toute mesure de probabilit\'e invariante 
par $f$ est inf\'erieure \`a $\delta$, alors il existe une 
fonction continue $u$ partout satisfaisant l'in\'egalit\'e 
\esp $|u - u \circ \Phi (f) - \log D\Phi (f) | \leq \varepsilon$
\esp pour tout $f \in \mathcal{G}$.}

\vspace{0.286cm}

La preuve r\'esulte  d'une application directe de la proposition suivante 
au cocycle $f \mapsto \log D \Phi (f)$ au dessus de l'action $\Phi$. 
(Rappelons qu'un {\em cocycle} associ\'e \`a une action d'un groupe 
par hom\'eomor- phismes d'un espace $X$ est une application 
$g \mapsto c(g) \in C(X)$ telle que l'\'egalit\'e $c(fg) = c(g) + c(f) \circ g$ 
est partout v\'erifi\'ee pour tout $f,g$ dans le groupe.)

\vspace{0.3cm}

\noindent{\bf Proposition.} {\em Soit $\Gamma$ un groupe nilpotent engendr\'e par 
une partie finie $\mathcal{G}$ et agissant par hom\'eomor- phismes sur un espace m\'etrique 
compact $X$. Il existe une partie g\'en\'eratrice finie $\mathcal{G}'$ contenant 
$\mathcal{G}$ telle que pour tout $\varepsilon > 0$ il existe $\delta > 0$ 
v\'erifiant la propri\'ete suivante : si  
$c: \Gamma \rightarrow C(X)$ est un cocycle associ\'e \`a cette action satisfaisant  
$\left| \int_X c(f) d \mu \right| < \delta$ pour tout $f \in \mathcal{G}'$ 
et toute probabilit\'e $\mu$ invariante par $f$, alors il existe $u \in C(X)$ 
telle que l'in\'egalit\'e $|u - u \circ f - c(f)| \leq \varepsilon$ est partout v\'erifi\'ee 
pour tout $f \in \mathcal{G}$.}

\vspace{0.21cm}

\noindent{\bf Preuve.} Tout d'abord, nous affirmons que si $c(f)$ a une moyenne 
de valeur absolue inf\'erieure \`a $\delta$ pour toute probabilit\'e invariante par $f$, 
alors il existe $N_0 \in \mathbb{N}$ tel que pour tout $N \geq N_0$ et tout $y \in X$,  
\begin{equation}\label{limit}
\left| \frac{1}{N} \sum_{n=0}^N c(f) \circ f^n (y) \right| < \delta.
\end{equation}
En effet, l'expression en question n'est autre que 
$$\int_X c(f) \hspace{0.04cm} d\mu_{y,N}, \qquad \mbox{o\`u} \qquad 
\mu_{y,N} := \frac{1}{N} \sum_{n=0}^{N-1} \delta_{f^n(y)}.$$
Puisque la famille $\mu_{y,N}$ contient des probabilit\'es invariantes par 
$f$ dans son adh\'erence, un argument simple de de compacit\'e permet 
de conclure que (\ref{limit}) a lieu pour tout $N$ assez large. 

Nous utiliserons le fait que les groupes nilpotents satisfont une 
version forte de la propri\'et\'e de {\em g\'en\'eration born\'ee} : il existe une partie g\'en\'eratrice 
finie $\mathcal{G}' = \{ f_1 , \ldots , f_{\ell} \}$ (qui peut \^etre choisie contenant $\mathcal{G}$) 
et une constante $M$ telles que tout \'el\'ement 
$f \in B(k) \subset \Gamma$ peut \^etre \'ecris sous la forme $f = f_{i_1}^{n_1} \cdots f_{i_m}^{n_m}$, 
o\`u chaque $f_{i_j}$ appartient \`a $\mathcal{G}'$ et de plus $m \leq M$ et $|n_j| \leq M k$  
(voir par exemple \cite[Appendix B]{BGr}). Dans ce qui suit, nous travaillerons avec la 
partie g\'en\'eratrice $\mathcal{G}'$, mais nous la noterons encore par $\mathcal{G}$.

Nous utiliserons aussi le fait que $\Gamma$ est \`a croissance polynomiale, ce qui implique 
l'existence d'une constante $C > 0$ ainsi que d'une suite croissante d'entiers positifs $k_n$ 
telles que 
$$\frac{\big| B(k_n+1) \setminus B(k_n) \big|}{\big| B(k_n) \big|} \leq \frac{C}{k_n},$$
o\`u $B(k)$ d\'esigne la boule de rayon $k$ dans 
$\Gamma$. En effet, dans le cas contraire, pour tout $d \geq 1$ on aurait pour 
quelques constantes positives $C',C'',C'''$ et tout $k$ suffissament large,
$$\big| B(k) \big| 
\geq C' \prod_{j=1}^k \Big(1 + \frac{C''}{j} \Big) 
\geq C''' \exp \Big( \sum_{j=1}^{k} \frac{C'''}{j} \Big)
\geq C''' k^d,$$
ce qui contredit la croissance polynomiale. 

Avec ces deux outils, nous pouvouns raisonner comme dans la cas commutatif au d\'ebut du 
\S III, avec quelques modifications. Posons 
$$u_n := \frac{1}{\big| B(k_n) \big|} \sum_{f \in B(k_n)} c(f).$$
\vspace{-0.1cm}Alors pour tout $f_i \in \mathcal{G}$ nous avons
\begin{eqnarray*}
u_n \circ f_i - u_n 
&=& \frac{1}{\big| B(k_n) \big|} \sum_{f \in B(k_n)} 
\big[ c(f) \circ f_i - c(f) \big]\\ 
&=& \frac{1}{\big| B(k_n) \big|} \sum_{f \in B(k_n)} 
\big[ c(f f_i) - c(f_i) - c(f) \big]
= -c(h_i) + \frac{1}{\big| B(k_n) \big|} \sum_{f \in B(k_n)} \big[ c(f f_i) - c(f) \big].
\end{eqnarray*}
Donc, la valeur de \hspace{0.03cm} 
$\big| u_n - u_n \circ f_i - c(f_i) \big|$ \hspace{0.03cm} est 
born\'ee par 
\begin{eqnarray*}
\frac{1}{\big| B(k_n) \big|}  \left| \sum_{f \in B(k_n)} \big[ c(f f_i) - c(f) \big] \right|
\leq \frac{1}{\big| B(k_n) \big|}  \sum_{f \in B(k_n+1) \setminus B(k_n)} \big| c(f) \big|
\end{eqnarray*}
En \'ecrivant chaque $f \in B(k_n+1) \setminus B(k_n)$ sous la forme 
$f = f_{i_1}^{n_1} \cdots f_{i_m}^{n_m}$ ci-dessus, cette expression 
se transforme en 
\begin{small}
$$\frac{1}{\big| B(k_n) \big|} \! \!\sum_{f \in B(k_n+1) \setminus B(k_n)} \!\!
\Big| \sum_{j=1}^{m} c(f_{i_j}^{n_j}) \circ f_{i_{j+1}}^{n_{j+1}} \cdots f_{i_n}^{n_m} \Big| 
=
\frac{1}{\big| B(k_n) \big|} \!\! \sum_{f \in B(k_n+1) \setminus B(k_n)} \!\!
\Big| \sum_{j=1}^{m} 
\sum_{n=0}^{n_j-1} c(f_{i_j}) \circ f_{i_j}^{n} f_{i_{j+1}}^{n_{j+1}} \cdots f_{i_n}^{n_m} \Big|. 
$$
\end{small}Puisque $m \leq M$, en notant $y:= f_{i_{j+1}}^{n_{j+1}} \cdots f_{i_n}^{n_m}(x)$ 
nous voyons que nous devons traiter des expressions du type 
$$\frac{1}{\big| B(k_n) \big|} \!\! \sum_{f \in B(k_n+1) \setminus B(k_n)} \!\!
\Big| \sum_{n=0}^{n_j-1} c(f_{i_j}) \circ f_{i_j}^{n} (y) \Big|.$$ 
Prenons $N_0$ de telle dorte que (\ref{limit}) ait lieu pour tout $N \geq N_0$, 
tout $y \in X$ et tout $f_i \in \mathcal{G}$. Deux cas peuvent alors se 
pr\'esenter. Si $n_j \leq N_0$ alors 
\begin{small}
$$\frac{1}{\big| B(k_n) \big|} \!\! \sum_{f \in B(k_n+1) \setminus B(k_n)} \!\!
\Big| \sum_{n=0}^{n_j-1} c(f_{i_j}) \circ f_{i_j}^{n} (y) \Big|
\leq 
\frac{N_0}{\big| B(k_n) \big|} \!\! \sum_{f \in B(k_n+1) \setminus B(k_n)} \!\! 
\max_{x\in X} \big| c(f_j)(x) \big| \leq \frac{C N_0 \max_{x\in X} \big| c(f_j)(x) \big|}{k_n},$$
\end{small}qui tend vers z\'ero lorsque $n$ tend vers l'infini. Si $n_j \geq N_0$ alors 
$$\frac{1}{\big| B(k_n) \big|} \!\! \sum_{f \in B(k_n+1) \setminus B(k_n)} \!\!
\Big| \sum_{n=0}^{n_j-1} c(f_{i_j}) \circ f_{i_j}^{n} (y) \Big|
< 
\frac{1}{\big| B(k_n) \big|} \!\! \sum_{f \in B(k_n+1) \setminus B(k_n)} \!\! \delta n_j 
\leq 
\frac{C\delta n_j}{k_n} \leq MC\delta,$$
qui est plus petit que $\varepsilon$ pour $\delta$ assez petit. 
$\hfill\square$

\vspace{0.362cm}

D'apr\`es ce qui pr\'ec\`ede, pour conclure la preuve du 
Th\'eor\`eme B dans le cas (i) o\`u il existe des orbites finies pour l'action, 
notre probl\`eme consiste \`a montrer que par conjugaison topologique on peut obtenir 
des conjugu\'es $C^1$ dont les moyennes (par rapport \`a toutes les probabilit\'es 
invariantes) du logarithme de la d\'eriv\'ee des g\'en\'erateurs sont petites en 
valeur absolue. Pour ce faire, notons d'une part qu'il suffit de consid\'erer les mesures
invariantes {\em ergodiques}, c'est-\`a-dire celles qui ne peuvent pas \^etre
exprim\'ees comme combinaison convexe non triviale de probabilit\'es invariantes. 
Puisque $\Gamma$ admet des orbites finies, tout \'el\'ement $f \! \in \! \Gamma$ poss\`ede 
des points p\'eriodiques. De plus, toute probabilit\'e invariante par $f$ doit \^etre support\'ee 
sur ses orbites p\'eriodiques, car le compl\'ementaire est form\'e par des points errants. 
En particulier, toute probabilit\'e invariante par $f$ et ergodique est la moyenne 
sur une telle orbite. L'int\'egrale de \hspace{0.04cm}
$\log (Df)$ \hspace{0.04cm} par rapport \`a une telle mesure n'est autre
que la moyenne des logarithmes des d\'eriv\'ees le long de l'orbite. 
Nous avons ainsi r\'eduit le probl\`eme \`a prouver le lemme suivant. 


\vspace{0.29cm}

\noindent{\bf Lemme (d'applatissement des points hyperboliques).} 
{\em Soit $\Gamma$ un sous-groupe de $\mathrm{Diff}_+^1(X)$ engendr\'e 
par une partie finie $\mathcal{G}$. 
Alors pour tout $\delta \!>\! 0$ il existe des conjugu\'es de $\Gamma$ 
par des hom\'eomorphismes $\psi$ tels que:  

\vspace{0.05cm}

\noindent -- pour tout $f \in \Gamma$, l'application \hspace{0.041cm}
$\psi \circ f \circ \psi^{-1}$ \hspace{0.041cm} est un diff\'eomorphisme 
de classe $C^1$, 

\vspace{0.05cm}

\noindent -- si $f \in \mathcal{G}$ admet un point p\'eriodique $x$ 
et $N$ est sa p\'eriode, 
alors le multiplicateur du conjugu\'e de $f$ au point $x$ satisfait 
\hspace{0.04cm} $e^{-\delta} \leq D(\psi \circ f^N \circ \psi^{-1})(x) \leq e^{\delta}$.}

\vspace{0.16cm}

\noindent{\bf Preuve.} Soit $\delta > 0$ donn\'e, et soit $\{x_1,\ldots,x_N\}$ 
l'une des orbites p\'eriodiques de $f \in \mathcal{G}$. Fixons $\alpha > 0$, et 
prenons $\psi_{\alpha} \in \mathrm{Hom\acute{e}o}_+(\mathrm{S}^1)$ qui soit 
un diff\'eomorphisme en dehors de $\{x_1,\ldots,x_N\}$ et pr\`es de chaque 
$x_j$ co\"{\i}ncide avec \hspace{0.02cm}
$y \mapsto x_j + (y-x_j)^{\frac{1}{\alpha}}.$ \hspace{0.02cm}
Nous affirmons alors que pour tout $g \in \Gamma$, chaque
$g_{\alpha} := \psi_{\alpha} \circ g \circ \psi_{\alpha}^{-1}$
est un diff\'eomorphisme de classe $C^1$ ; de plus, si $\alpha$ 
est suffisamment grand, alors pour tout indice $j$,
$$\big| \log D(\psi_{\alpha} \circ f^N \circ \psi_{\alpha}^{-1}) (x_j) \big| \leq \delta.$$
En effet, on constate d'abord que pour $y$ proche de $x_j$,
$$\psi_{\alpha} \circ g \circ \psi_{\alpha}^{-1} (y)
= \big[ g \big( x_j + (y-x_j)^{\alpha} \big) - g(x_j) \big]^{1/\alpha} + g (x_j).$$
Donc,
\vspace{-0.4cm}
\begin{eqnarray*}
D (\psi_{\alpha} \!\circ\! g \!\circ\! \psi_{\alpha}^{-1}) (y)
\!\!\!&=& \!\!\!\alpha (y-x_j)^{\alpha - 1} Dg \big (x_j + (y-x_j)^{\alpha} \big)
\frac{1}{\alpha} \big[g ( x_j + (y-x_j)^{\alpha}) - g (x_j) \big]^{\frac{1}{\alpha} - 1}\\
&=& \!\!\! (y-x_j)^{\alpha - 1} Dg \big (x_j + (y-x_j)^{\alpha} \big)\!\!
\left[ \frac{g ( x_j + (y-x_j)^{\alpha}) - g(x_j)}{(y-x_j)^{\alpha}} \right]^{\frac{1}{\alpha} - 1}
\!\! \big[ (y-x_j)^{\alpha} \big]^{\frac{1}{\alpha} - 1}\\
&=& \!\!\! Dg \big (x_j + (y-x_j)^{\alpha} \big)
\left[ \frac{g ( x_j + (y-x_j)^{\alpha}) - g(x_j)}{(y-x_j)^{\alpha}} \right]^{\frac{1}{\alpha} - 1}.
\end{eqnarray*}
Lorsque $y$ tend vers $x_j$, cette expression converge vers
$$Dg(x_j) \big[ Dg(x_j) \big]^{\frac{1}{\alpha}-1} = \big[ Dg(x_j) \big]^{\frac{1}{\alpha}},$$
ce qui d\'emontre de mani\`ere simultan\'ee les deux propri\'et\'es annonc\'ees.

Nous pouvouns r\'ep\'eter cet argument avec chaque orbite p\'eriodique 
d'un g\'en\'erateur de $\Gamma$ dont le multiplicateur soit
$\leq e^{-\delta}$ ou $\geq e^{\delta}$. Puisqu'il n'y a qu'un 
nombre fini de telles orbites, ceci permet de conclure la preuve du lemme. 
$\hfill\square$

\vspace{0.32cm}

La preuve du Th\'eor\`eme B est enfin termin\'ee. 

\vspace{0.2cm}

\noindent{\bf Remarque.} Ci-dessus, on aurait pu conjuguer directement par un
hom\'eomorphisme qui est un diff\'eomorphisme loin des points de l'orbite concern\'ee
et dont le germe autour de ces points co\"{\i}ncide avec celui de $x \to \exp(-1/x)$
\`a l'origine. En effet, ceci rend tangent \`a l'identit\'e tout diff\'eomorphisme
de classe $C^1$ qui pr\'eserve l'orbite (voir \cite{tsuboi}). 

\vspace{0.3cm}


\noindent{\bf IV. Sur la connexit\'e par arcs.} Dans la preuve du Th\'eor\`eme B,
l'action par des translations est aboutie par un chemin continu de r\'epr\'esentations
lorsqu'il n'y a pas d'orbite finie. Ceci est encore valable lorsqu'il y a des orbites
finies mais les points p\'eriodiques des g\'en\'erateurs sont tous paraboliques. 
En effet, la preuve du lemme du \S II. s'applique encore dans ce contexte.

Dans le cas qui nous reste, c'est-\`a-dire lorsque quelques g\'en\'erareurs 
poss\`edent des orbites p\'eriodiques hyperboliques, nous voudrions 
encore aboutir \`a une r\'epr\'esentation par des translations par un 
chemin form\'e par des conjugu\'es topologiques de l'action originelle. Pour
ce faire, on aimerait appliquer une m\'ethode semblable \`a celle de la fin du \S III. en
prenant une famille de diff\'eomorphimes $\psi_{\alpha}$ qui varie contin\^ument
pour la topologie $C^1$ par rapport \`a $\alpha$ en dehors de telles orbites (ou
du moins, en dehors des orbites avec un multiplicateur trop grand). En effet, les
calculs de la preuve du lemme d'applatissemment des points hyperboliques 
montrent que pour tout $f \in \Gamma$, l'application
$\alpha \to f_{\alpha}$ est continue pour la topologie $C^1$. Cependant, ce
proc\'ed\'e pourrait {\em \`a priori} faire exploser les d\'eriv\'ees. Nous
devrons donc \^etre un peu plus soigneux, et pour cela il nous sera plus confortable
de raisoner au niveau des \'equations cohomologiques, tout en cherchant des solutions
rapproch\'ees non n\'ecessairement born\'ees au vosinage des extr\'emit\'es.

\vspace{0.3cm}

\noindent{\bf Fin de la preuve du Th\'eor\`eme C.} 
Soit $\Gamma$ un groupe nilpotent engendr\'e par une partie finie $\mathcal{G}$ 
et admettant des orbites finies. \'Etant donn\'e $\varepsilon > 0$, fixons le 
$\delta > 0$ donn\'e par le lemme d'existence de solutions \'erron\'ees. 
Bien \'evidemment, la r\'eunion des orbites des
points p\'eriodiques hyperboliques des g\'en\'erateurs \`a multiplicateur
soit $\leq e^{-\delta}$ soit $\geq e^{\delta}$ est un ensemble fini, disons
$\{x_1,\ldots,x_k\}$. Soit $\psi$ un hom\'eomorphisme de $X$ qui est un
diff\'eomorphisme $C^1$ restreint au compl\'ementaire de $\{x_1,\ldots,x_k\}$ 
et fixe chaque point $x_i$ de telle sorte que son germe autour d'un tel point
co\"{\i}ncide avec celui de $x \to x^{\alpha}$ \`a l'origine. Pour $\alpha > 0$
assez large, les multiplicateurs sur les points $x_i$ des conjugu\'es par $\psi$ des
g\'en\'erateurs sont tous compris entre $e^{-\delta}$ et $e^{\delta}$. Fixons
un tel $\alpha$ et notons $v : = \log (D \psi)$. Puisque l'action conjugu\'ee
par $\psi$ est une action par diff\'eomorphismes de classe $C^1$ qui satisfait 
les hypoth\`eses du lemme d'existence de solutions \'erronn\'ees, il existe une 
solution $\varepsilon$-rapproch\'ee continue $w$ de l'\'equation cohomologique
associ\'ee \`a cette action, que l'on peut prendre comme \'etant le
logarithme de la d\'eriv\'ee d'un diff\'eomorphisme $\varphi$ de classe
$C^1$ de $X$. La fonction $u : = w \circ \psi + v$ est alors une solution
$\varepsilon$-raproch\'ee de l'\'equation cohomologique associ\'ee \`a l'action
originelle. M\^eme si elle n'est pas continue (elle est non born\'ee au voisinage
des points $x_i$), elle co\"incide sur $X \setminus \{x_1,\ldots,x_k\}$ avec
le logarithme de la d\'eriv\'ee de $\phi := \varphi \circ \psi$, qui est
hom\'eomorphisme de $X$ qui transforme par conjugaison l'action originelle
en une autre action par diff\'eomorphismes de classe $C^1$ de $X$.

Pour chaque $\varepsilon := 1/n$ prenons la fonction $u := u_n$ induite
comme ci-dessus. Ces fonctions sont reli\'ees par des chemins affines
\esp $t \to (1-t) u_n + t u_{n+1} =: \hat{u}_{n+t}$. \esp Nous affirmons
que quitte \`a rajouter une constante $C_{n+t}$, on peut supposer que
$u_{n+t} := \hat{u}_{n+t} + C_{n+t}$ est le logarithme de la
d\'eriv\'ee d'un hom\'eomorphisme $\phi_{n+t}$ de $X$.
En effet, l'int\'egrale totale de la fonction $\exp (\hat{u}_{n+t})$ est
finie : ceci d\'ecoule directement de la convexit\'e de la fonction exponentielle
en tenant compte que $\exp(u_n)$ et $\exp(u_{n+1})$ sont toutes les deux
d'int\'egrale totale finie (\'egale \`a 1). On prend alors $C_{n+t}$ comme
\'etant l'exponentielle de l'inverse de cette int\'egrale, on fixe
$x_0 \in X$ et on d\'efinit
$$\phi_{n+t}(x) := \int^{x}_{x_0} \exp u_{n+t}.$$
Le logarithme de la d\'eriv\'ee du conjugu\'e par $\phi_{n+t}$ de
$f \in \Gamma$ au point $\phi^{-1} (x)$ est \'egal \`a
\begin{small}
$$\log Df (x) + u_{n+t} \circ f (x) - u_{n+t} (x) =
(1-t) \big[ \log Df (x) + u_{n} \circ f (x) - u_n (x) \big]
+ t \big[ \log Df (x) + u_{n+1} \circ f (x) - u_{n+1} (x) \big].$$
\end{small}Il s'agit donc une fonction continue, ce qui montre que l'action originelle
conjugu\'ee par $\phi_{n+t}$ est une action par diff\'eomorphismes de classe $C^1$.
De plus, cette action varie contin\^ument par rapport au param\`etre $t$. Finalement,
l'\'egalit\'e ci-dessus montre que pour tout $f \in \mathcal{G}$ le logarithme de la
d\'eriv\'ee de $\phi_{n+t} \circ f \circ \phi_{n+t}^{-1}$ est inf\'erieur ou \'egal
\`a $1/n$ partout. Il s'en suit que le chemin d'actions conjugu\'ees par $\phi_{n+t}$
est continu pour la topologie $C^1$ et aboutit \`a l'infini en une r\'epr\'esentation
par des translations, tel qu'on le d\'esirait.


Nous venons donc de montrer que toute r\'epr\'esentation soit dans 
$\mathrm{Diff}_+^1([0,1])$ soit dans $\mathrm{Diff}_+^1 (\mathrm{S}^1)$ 
est dans la m\^eme composante connexe par arcs d'une r\'epr\'esentation par 
des translations. En particulier, \'etant donn\'ees deux r\'epr\'esentations dans 
$\mathrm{Diff}^1_+ ([0,1])$, on peut les joindre par un chemin qui passe par la 
r\'epr\'esentation triviale. Ceci montre la connexit\'e par arcs pour le cas des actions 
sur l'intervalle. Pour le cas du cercle, on r\'ep\`ete l'argument de la fin du \S II. 
$\hfill\square$


\vspace{0.53cm}

\noindent{\bf V. Quelques commentaires finales.} 
Nous voudrions conclure par une proposition et un corollaire qui aident \`a mieux 
comprendre le cas des orbites finies. M\^eme si l'on peut donner des preuves purement 
combinatoires de ces deux r\'esultats, mais nous pr\'ef\'erons de faire appel \`a 
\cite{ollagnier-pinchon} afin d'illustrer le type d'information contenue dans le 
r\'esultat cohomologique l\`a-dedans.

\vspace{0.2cm}

\noindent{\bf Proposition (\`a propos des points p\'eriodiques hyperboliques).}
{\em Si $\Gamma$ est un groupe nilpotent de type fini de diff\'eomorphismes de
classe $C^1$ de $X$, alors les deux conditions suivantes sont \'equivalentes :}

\noindent -- {\em les points p\'eriodiques des \'el\'ements de $\Gamma$
sont tous paraboliques,}

\noindent -- {\em les points p\'eriodiques des \'el\'ements de $\Gamma$
appartenant au support d'une mesure de probabilit\'e invariante par
$\Gamma$ sont tous paraboliques.}

\vspace{0.1cm}

\noindent{\bf Preuve.} Si la deuxi\`eme condition est valable, alors d'apr\`es 
\cite[Th\'eor\`eme 1]{ollagnier-pinchon} nous sommes encore dans
le cas o\`u il existe des solutions rapproch\'ees \`a notre \'equation cohomologique.
Or, l'existence de telles solutions entra\^{\i}ne que les points p\'eriodiques
des \'el\'ements sont tous paraboliques. En effet, une in\'egalit\'e (partout)
du type
$$\big| u - u \circ f - \log (Df) \big| < \varepsilon$$
entra\^{\i}ne
$$\big| u - u \circ f^N - \log (D f^N) \big| < N \varepsilon.$$
En particulier, si $f^N (x_0) = x_0$, alors \hspace{0.04cm} 
$\big| \log D f^N (x_0) \big| < N \varepsilon.$ \hspace{0.04cm}
Donc, si de telles fonctions $u$ existent pour tout $\varepsilon > 0$,
alors nous avons n\'ecessairement  $Df^N (x_0) = 1$. $\hfill\square$

\vspace{0.3cm}

\noindent{\bf Corollaire.} {\em Soit $\Gamma$ un groupe nilpotent de type fini
de diff\'eomorphismes de classe $C^1$ de $X$. Si $f \in \Gamma$ admet un point
p\'eriodique hyperbolique $x_0$, alors l'orbite de $x_0$ par $\Gamma$ est finie.}

\vspace{0.1cm}

\noindent{\bf Preuve.} D'apr\`es le lemme du \S II. et sa d\'emonstration, $\Gamma$ admet
des orbites finies ; de plus, le sous-ensemble $\hat{\Gamma}$ des \'el\'ements de $\Gamma$
ayant des points fixes est un sous-groupe distingu\'e et d'indice fini (c'est le noyau de la 
fonction nombre de translation, qui est dans ce cas un homomorphisme ; 
voir \cite[Lemma 4.1.2]{navas-book}). Si $N$ est la
p\'eriode de $x_0$ pour $f$, alors $f^N$ appartient \`a $\hat{\Gamma}$. Supposons
que $x_0$ ne soit pas fix\'e par $\hat{\Gamma}$, et soient $a,b$ les points fixes de
$\hat{\Gamma}$ \`a gauche et \`a droite de $x_0$ respectivement. 
Alors la proposition pr\'ec\'edente
donne une contradiction lorsqu'on l'applique \`a l'action de $\hat{\Gamma}$ sur $[a,b]$
apr\`es conjugaison par un hom\'eomorphisme qui applatit compl\'etement les d\'eriv\'ees
aux extr\'emit\'es (voir la remarque du \S III.). Par suite, $\hat{\Gamma}$ fixe $x_0$, et
puisque $\Gamma / \hat{\Gamma}$ est fini cyclique, l'orbite de $x_0$ par $\Gamma$
doit \^etre finie. $\hfill\square$

\vspace{0.5cm}

\noindent{\bf Remerciements.} Je remercie S.~Crovisier, \'E.~Ghys, E.~Giroux et
A.~Kocsard, qui ont \'et\'e les premiers \`a me signaler la question de la connexit\'e
des r\'epr\'esentations de $\mathbb{Z}^d$ par diff\'eomorphismes en dimension 1. 
Je remercie aussi  C.~Bonatti, H.~Eynard, N.~Guelman  et \'E.~Farinelli, dont les
r\'esultats gentiment communiqu\'es ont \'et\'e une source importante
d'inspiration, ainsi que Y. Matsuda, qui m'a signal\'e une erreur dans 
l'\'enonc\'e du Th\'eor\`eme C de la premi\`ere version de ce travail. 
Finalement, je remercie le rapporteur anonyme dont les nombreuses 
remarques et critiques ont permis de nettement am\'eliorer la pr\'esentation.

Cette Note a \'et\'e \'ecrite lors de s\'ejours \`a
l'Institut des Math\'ematiques Pures de T\'eh\'eran et \`a l'Institut Mittag
Leffler. Je remercie M.~Nassiri et A.~Karlsson pour leurs invitations respectives.

Ce travail a \'et\'e financ\'e par le projet ACT 1103 DySyRF (Center of Dynamical 
Systems and Related Fields) ainsi que par le projet FONDECYT 1120131.


\begin{small}



\vspace{0.3cm}

\noindent Andr\'es Navas

\noindent Dpto. de Matem\'atica y C.C., Univ. de Santiago de Chile (USACH)

\noindent Alameda 3363, Estaci\'on Central, Santiago, Chile

\noindent Email: andres.navas@usach.cl

\end{small}

\end{document}